\documentclass[12pt,a4paper]{article}
\usepackage{amsmath,amsthm,amssymb,latexsym,epic}
\usepackage{graphicx}

\newtheorem{theorem}{Theorem}
\newtheorem{lemma}{Lemma}
\newtheorem{corollary}{Corollary}
\newtheorem{proposition}{Proposition}

\usepackage[all]{xy}
\usepackage{enumerate}
\newcommand{\N}{{\mathbb N}}

\newcommand{\dom}{\mathrm{dom}}
\renewcommand{\ker}{\mathrm{ker}}
\newcommand{\ran}{\mathrm{ran}}
\renewcommand{\L}{\mathcal{L}}
\newcommand{\R}{{\mathcal R}}
\newcommand{\ISA}{{\mathcal {ISA}}}
\renewcommand{\H}{{\mathcal H}}
\newcommand{\D}{{\mathcal D}}
\newcommand{\J}{{\mathcal J}}
\newcommand{\T}{{\mathcal T}}
\newcommand{\PT}{{\mathcal {PT}}}
\newcommand{\IS}{{\mathcal {IS}}}
\renewcommand{\S}{{\mathcal {S}}}

\newcommand{\PAut}{\mathrm{PAut}}
\newcommand{\PEnd}{\mathrm{PEnd}}
\newcommand{\End}{\mathrm{End}}

\begin{document}
\markboth{G. Kudryavtseva} {On conjugacy in regular epigroups}

%
%

\title{On conjugacy in regular epigroups}
\author{Ganna Kudryavtseva}
\date{}
\maketitle
\begin{center}{Department of Mechanics and Mathematics, \\
Kyiv Taras Shevchenko University,\\ Volodymyrs'ka, 60, Kuiv, 01033, Ukraine\\
akudr@univ.kiev.ua}
\end{center}

\begin{abstract}
Let $S$ be a semigroup. The elements $a,b\in S$ are called
primarily conjugate if $a=xy$ and $b=yx$ for certain $x,y\in S$.
The relation of conjugacy is defined as the transitive closure of
the relation of primary conjugacy. In the case when $S$ is a
monoid, denote by $G$ the group of units of $S$. Then the relation
of $G$-conjugacy is defined by $a\sim_G b \iff a=g^{-1}bg$ for
certain $g\in G$. We establish the structure of conjugacy classes
for regular epigroups (i.e. semigroups such that some power of
each element lies in a subgroup). As a corollary we obtain a
criterion of conjugacy in terms of $G$-conjugacy for factorizable
inverse epigroups. We show that our general conjugacy criteria
easily lead to known and new conjugacy criteria for some specific
semigroups, among which are the full transformation semigroup and
the full inverse symmetric semigroup over a finite set, the linear
analogues of these semigroups and the semigroup of finitary
partial automatic transformations over a finite alphabet.
\end{abstract}

AMS 2000 Subject Classification 20M10, 20M20.
\date{}
\maketitle

\section{Introduction}\label{s1}

The notion of conjugacy in semigroups can be generalized from the
corresponding notion for groups in several ways. Perhaps, the  two
most natural and commonly used notions are the relations $\sim_G$
and $\sim$, whose definitions below are taken from~\cite{La}.

Let $S$ be a monoid and $G$ the group of units of $S$. The
relation $\sim_G$, called {\em $G$-conjugacy}, is defined as
$a\sim_G b$ if and only if $a=g^{-1}bg$ for certain $g\in G$. Let
now $S$ be a semigroup. We call the elements $a,b\in S$ {\em
primarily $S$-conjugate} if there exist $x,y\in S$ such that
$a=xy,$ $b=yx$. This will be denoted by $\sim_{pS}$ or just by
$\sim_p$ when this does not lead to ambiguity. The relation
$\sim_p$ is reflexive and symmetric while not transitive in the
general case. Denote by $\sim$ the transitive closure of the
relation $\sim_p$. If $a\sim b$ then $a$ and $b$ are said to be
{\em $S$-conjugate} or just {\em conjugate}. It is easy to see
that in the case of group, both $\sim_G$ and $\sim$ coincide with
usual group conjugacy. Besides, for a monoid $S$ there is an
inclusion $\sim_G\subset \sim$.

The structure of conjugacy and $G$-conjugacy classes for some
specific regular semigroups was studied in a number of papers
(\cite{GK}, \cite{KM}, \cite{KM1}, \cite{Ch}, \cite{OS}), see also
the monograph \cite{Li}. For the relation of conjugacy the
structure of conjugacy classes usually happens to be more
complicated than for the relation of $G$-conjugacy.

The present paper is devoted to the systematic study of the
conjugacy relation in regular epigroups (an {\em epigroup} or a
{\em group-bound semigroup} is a semigroup such that some power of
each its element lies in a subgroup) and is organized as follows.
In the Preliminaries we collect some notation used throughout the
paper and cite some well-known facts about the structure of
$\D$-classes. In Section~\ref{s3} we establish a criterion of
conjugacy of two group elements of a given semigroup. In
section~\ref{s4} we establish Theorem~\ref{main}, which gives a
criterion of conjugacy of two group-bound elements of a regular
semigroup and Theorem~\ref{p3}, which provides a criterion of
conjugacy in terms of $G$-conjugacy for factorizable inverse
epigroups. In Section~\ref{appl} we show that conjugacy criteria
for many important specific examples of regular epigroups can be
derived in a unified way from our main results. In particular, we
give short and very clear proofs for some known conjugacy
criteria. Besides,  for the first time we formulate and prove
conjugacy criteria for the semigroups $\PAut(V)$, $\PEnd(V)$,
$Fin\ISA(X)$, $Fin{\mathcal A}(X)$ and $Fin{\mathcal {PA}}(X)$
(see the explanations for the notations at the appropriate places
of the paper).

In the case when a regular semigroup $S$ is not group-bound, the
problem of description of conjugacy classes of $S$ seems to be
much more complicated, in particular, even for such a classical
semigroup as $\T({\mathbb N})$ the conjugacy classes are not
classified yet. At the same time, for the semigroup $\IS({\mathbb
N})$, which is also not an epigroup, the conjugacy classes are
described (see~\cite{KM}).

In Appendix~A  we show that despite finiteness of $X$ the
semigroup $\ISA(X)$ of partial automatic permutations over a
finite alphabet $X$ ($|X|\geq 2$) is not an epigroup. This
implies, in particular, that in $\ISA(X)$ there are conjugacy
classes without group elements showing that the conjugacy
criterion for the semigroup $\ISA(X)$ announced in Theorems 3,4
of~\cite{OS} fails to give sufficient condition of conjugacy. At
the same time we have substantial arguments that the description
of conjugacy classes in $\ISA(X)$ can be obtained using the
methods from~\cite{GNS, KM}. A paper devoted to this question is
now in preparation.

\section{Preliminaries}

Let $S$ be a semigroup and $a\in S$. The class containing $a$ with
respect to the  $\H$- ($\L$-, $\R$-, $\D$-, $\J$-) Green's
relations will be denoted by $H_a$ ($L_a$, $R_a$, $D_a$, $J_a$).
The following well-known facts about the structure of $\D$-classes
will be used and referred to in the sequel. Their proofs can be
found, for example, in~\cite{Hig}.

\begin{proposition}[see \cite{Hig}, Theorem 1.2.5, p. 18]\label{d1} Let $a,b\in S$. Then $ab\in R_a\cap L_b$ if and
only if $R_b\cap L_a$ contains an idempotent. In particular, a
triple  $a,b, ab$ belongs to the same $\H-$ class if and only if
this $\H$-class is a group.
\end{proposition}

\begin{proposition}[see \cite{Hig}, Theorems 1.2.7 and 1.2.8., pp. 18, 19]\label{d3} Let $e,f \in S$ be idempotents and $e\D f$. Then
for any $t\in R_e\cap L_f$ there is an inverse $t'$ of $t$ such
that $t'\in R_f\cap L_e$. Furthermore, the maps $\rho_t\circ
\lambda_{t'}:$ $H_e\to H_f$ and $\rho_{t'}\circ \lambda_t:$
$H_f\to H_e$ defined via $x\mapsto t'xt$ and $x\mapsto txt'$,
respectively, are mutually inverse isomorphisms.
\end{proposition}

Recall that an element $a\in S$ is said to be a {\em group
element} provided $a$ belongs to a certain subgroup of $S$. It is
easily seen and well-known that for a group element $a\in S$ its
$\H$-class $H_a$ is a group (in fact, $H_a$ is a maximal subgroup
of $S$). Denote by $a^{-1}$ the (group) inverse of $a$ in $H_a$.

Let $a\in S$ and there exists $t\in\N$ such that $a^t$ is a group
element. In this case $a$ is called a {\em group-bound element}.
$S$ is called an {\em epigroup} (or a {\em group-bound semigroup})
provided that each element of $S$ is group-bound.

The following fact is known and is easily proved.

\begin{lemma}\label{ll1} The following statements are equivalent.
\begin{enumerate}
\item\label{one}  $a^k\H a^t$ for some $k>t$.
\item\label{two}  $a^i\H a^t$ for all $i\geq t$.
\item\label{three} $H_{a^t}$ is a group.
\end{enumerate}
\end{lemma}

Let $a\in S$ be a group-bound element and $t\in \N$ is such that
$H_{a^t}$ is a group. It follows from Lemma~\ref{ll1} that we can
correctly define $e_a$  (the notation goes from  \cite{Sh1}) to be
the identity element of the group $H_{a^t}$. Using Lemma~\ref{ll1}
one can easily obtain the following (known) useful statement.

\begin{corollary}\label{c1}
Suppose $a$ is a group-bound element of $S$. Then $e_aa=ae_a$ and
$ae_a\H e_a$. In particular, $ae_a$ is a group element.
\end{corollary}

\section{Conjugacy criterion for group elements}\label{s3}

We start from conjugacy criterion for group elements of an
arbitrary semigroup $S$ generalizing a similar result obtained
earlier in~\cite{KM} for the case when the semigroup $S$ is
finite.

Recall that elements $a,b\in S$ are said to be {\em mutually
inverse} provided that $a=aba$ and $b=bab$.

\begin{theorem}\label{pp1}
Let $S$ be a semigroup and $a,b\in S$ group elements. Then
\begin{enumerate}
\item $a\sim_p b$ if and only if there exists a pair of mutually inverse
elements $u,v\in S$ such that $b=uav$ and $a=vbu$.
\item  $a\sim b$ if and only if $a\sim_p b$.
\end{enumerate}
\end{theorem}

To prove this theorem we will need the following five lemmas.

\begin{lemma}\label{aux}
$a\sim_p b$ implies $a^n\sim_p b^n$ for each $n\geq 1$.
\end{lemma}

\begin{proof}
It is enough to note that if $a=xy$ and $b=yx$ then
$a^n=x(yx)^{n-1}\cdot y$ and $b=y\cdot x(yx)^{n-1}$, $n\geq 2$.
\end{proof}

\begin{lemma}\label{lem:new}
Suppose $a, b\in S$ and $a\sim_p b$. If $b$ belongs to a group,
then so does $a^2$.
\end{lemma}

\begin{proof}
Assume $a=ts$, $b=st$ for certain $t,s\in S$.  That
$b=e_{b}de_{b}=e_{b}ste_{b}$ implies $b\in e_{b}sS^1$. Besides,
$e_{b}s=bb^{-1}e_{b}s\in bS^1$, whence $b\R e_{b}s$. Analogously
one shows that $b\L te_{b}$. Therefore, $e_{b}s\cdot te_{b} \in
R_{e_{b}s}\cap L_{te_{b}}$ and thus in view of
Proposition~\ref{d1}
 $L_{e_{b}s}\cap R_{te_{b}}$ contains an
idempotent. Now, since $L_{te_{b}}\cap R_{e_{b}s}=H_{b}$ is a
group, it follows that $te_{b}\cdot e_{b}s = te_{b}s\in
R_{te_{b}}\cap L_{e_{b}s}$, so that $te_{b}s$ is a group element.
This implies $(te_{b}s)^2\H te_{b}s$. Therefore, since

$$a^{2}=tsts=t{b}s=te_{b}de_{b}s=(te_{b}s)^2,$$

we have that $a^{2}$ is a group element.
\end{proof}

Say that $c, d\in S$ are {\em conjugate in at most $k$ steps}
provided that there are $k\geq 1$ and $c=c_0, c_1\dots c_k=d$ such
that $c_i\sim_p c_{i+1}$, $0\leq i\leq k-1$.

\begin{lemma}\label{lem:new1}
Suppose $a, b\in S$ and $a\sim b$. If $b$ belongs to a group, then
so does some power of $a$.
\end{lemma}

\begin{proof}
Let $k\geq 1$ be such that $a$ and $b$ are conjugate in at most
$k$ steps. Show that $a^{2^k}$ is a group element. Apply induction
on $k$. For $k=1$ the statement follows from Lemma~\ref{lem:new}.

Assume now that $m\geq 1$ and the statement is proved for $k=m$.
Let $k=m+1$.  Fix $a=c_0, c_1,\dots, c_k=b$ such that $c_i\sim_p
c_{i+1}$, $0\leq i\leq k-1$. Since $c_{k-1}\sim_p b$ then
$c_{k-1}^2$ is a group element by Lemma~\ref{lem:new}. Beside
this, $a^{2}\sim_p c_1^2\sim_p\dots \sim_p c_{k-1}^2$ due to
Lemma~\ref{aux}, which means that $a^{2}$ and $c_{k-1}^2$ are
conjugate in at most $k-1=m$ steps. Applying the inductive
hypothesis we obtain that $(a^{2})^{2^{k-1}}=a^{2^k}$ is a group
element, as required.
\end{proof}

\begin{corollary}\label{epi}
Suppose $a, b\in S$ and $a\sim b$. If $b$ is group-bound then so
is $a$.
\end{corollary}

\begin{proof} The statement follows from Lemma~\ref{aux} and
Lemma~\ref{lem:new1}.
\end{proof}

\begin{lemma}\label{l1}
Suppose that $a,b\in S$ are group elements and $a\sim_p b$. Then
there exist $x,y\in S$ such that $a=xy$, $b=yx$ and $x\in R_a\cap
L_b$, $y\in R_b\cap L_a$.
\end{lemma}
\begin{proof}
Since $a\sim_p b$ then $a=st$ and $b=ts$ for certain $s,t\in S$.
It follows that $te_a\L a\R e_as$. Then $e_as\cdot te_a \in
R_{e_as}\cap L_{te_a}$ which implies that
 $L_{e_as}\cap R_{te_a}$ contains an
idempotent by Proposition~\ref{d1}. Since $L_{te_a}\cap
R_{e_as}=H_{a}$ is a group then $te_a\cdot e_as = te_as\in
R_{te_a}\cap L_{e_as}$, so that $te_as$ is a group element. This
implies  $(te_as)^2\H te_as$. Therefore, in view of
$$(te_as)^2=t\cdot
e_aste_a\cdot s=tsts=b^2\H b,$$ we get $te_as\H b$. Hence $te_a\R
te_as\H e_b$, whence $e_bte_a=te_a$. Analogously, $se_bt\H a$ and
$e_ase_b=e_as$.  But then
$$a=e_as\cdot te_a=e_ase_b\cdot e_bte_a=e_a\cdot se_bt\cdot
e_a=se_bt$$ and analogously $b=te_as$. Set $x=e_ase_b$,
$y=e_bte_a$. We obtain $a=xy$, $b=yx$ and $x\in R_a\cap L_b$,
$y\in R_b\cap L_a$ as required.
\end{proof}

\begin{lemma}\label{la1}
Let $a,b\in S$ be two group-bound elements. Then $a\sim_p b$
implies $ae_a\sim_p be_b$.
\end{lemma}
\begin{proof}
Let $n\in \N$ be chosen such that $a^n$ and $b^n$ are group
elements. Fix $x,y\in S$ such that $a=xy$ and $b=yx$. Then
$$
a^{n+1}=x(yx)^ny=xb^ny=xb^ne_by=a^nxe_by=a^ne_axe_by.
$$
Multiplying both sides of this equality  by $(a^n)^{-1}$ from the
left we obtain
$$e_aa=(a^n)^{-1}a^{n+1}=e_axe_by.$$
Since $e_aa=ae_a$ by Corollary~\ref{c1} it follows that
$ae_a=e_axe_by$. Similarly, $be_b=e_bye_ax$. Therefore,
$ae_a\sim_p be_b$.
\end{proof}

\begin{lemma}\label{l2} Let $a,b \in S$ be group elements satisfying $a\H b$ and $a\sim_p
b$. Then there exists $h\in H_a$ such that $a=h^{-1}bh$.
\end{lemma}
\begin{proof} By Lemma~\ref{l1} $a=hg$, $b=gh$ for some $h,g\in
H_a$. Thus, $ah^{-1}=h^{-1}b$, which implies $a=h^{-1}bh$ as
required.
\end{proof}

Now we are ready to prove Theorem~\ref{pp1}.

\begin{proof}[Proof of Theorem~\ref{pp1}]
{\em (1).} {\em Necessity.} Let $a\sim_p b$. Fix a pair of
mutually inverse elements $t\in R_a\cap L_b$ and $t'\in R_b\cap
L_a$ (this is possible to do by Proposition~\ref{d3}). In
particular, $e_a=tt'$, $e_b=t't$ . Fix also some $x\in R_a\cap
L_b$ and $y\in R_b\cap L_a$ such that $a=xy$ and $b=yx$ (such
elements exist by Lemma~\ref{l1}). Then $b\H t'at=t'x\cdot yt$ and
$yt\cdot t'x=ye_bx=yx=b$. It follows that $t'at\sim_p b$ and
$t'at\H b$. Now Lemma~\ref{l2} ensures us that there is $g\in H_b$
such that $b=g^{-1}t'atg$. Set $u=g^{-1}t'$, $v=tg$. Since
$L_g\cap R_{t'}=H_g$ contains an idempotent then $u=g^{-1}t'\in
R_g\cap L_{t'} =H_{t'}$. Similarly, $v\in H_{t}$. Furthermore,
$uvu=g^{-1}t'tgg^{-1}t'= u$, $vuv=tgg^{-1}t'tg=v$. Thereby, using
Proposition~\ref{d3}, we have that $\rho_u\circ\lambda_{v}$ is an
isomorphism from $H_a$ to $H_b$. It remains to note, that $b=uav$
and $a=vbu$.

{\em Sufficiency.} Suppose $b=uav$ and $a=vbu$, where $u,v$ are
mutually inverse. Then $b=uvbuv$ which implies that $uvb=uvbuv$.
The two previous equalities imply $b=uvbuv=uvb$. Denote $s=vb$,
$t=u$. Then $st=a$ and $ts=uvb=b$. Hence, $a\sim_p b$.

{\em (2).} Clearly, we have to show only that $a\sim b$ implies
$a\sim_p b$.  Suppose that $a$ and $b$ are conjugate in at most
$n$ steps and that $a=a_0, a_1,\dots, a_n=b$ such that $a_i\sim_p
a_{i+1}$, $0\leq i\leq n-1$ are fixed. Corollary~\ref{epi} implies
that all $a_i$ are group-bound. Apply induction on $n$. If $n=1$
there is nothing to prove.

Let $n=2$. Suppose $a\sim_p a_1$, $a_1\sim_p b$. By the first
statement of this Theorem $a=ta_1s$, where $ts=e_a$, $st=e_{a_1}$,
and $a_1=ubv$, where $uv=e_{a_1}$, $vu=e_b$. Then $a=tubvs$,
$b=va_1u=vsatu$ and $tuvs=te_{a_1}s=ts=e_a$,
$vstu=ve_{a_1}u=vu=e_b$.

Let $n\geq 2$. Assume that any two group elements, which are
conjugate in at most $k$ steps with $k\leq n-1$, are primarily
conjugate.  It follows from Lemma~\ref{la1} that
$$a=a_0e_{a_0}\sim_p a_1e_{a_1}\sim_p\dots\sim_p a_ne_{a_n}=b.$$

Note that all $a_ie_{a_i}$ are group elements by
Corollary~\ref{c1}.  Then $a\sim_pae_{a_{n-1}}$ by the inductive
assumption. It follows that $a\sim_pa_{n-1}e_{a_{n-1}}\sim_p b$.
The inductive assumption implies now that $a\sim_p b$.
\end{proof}

The following statements are direct consequences of
Theorem~\ref{pp1}.

\begin{corollary}\label{c2} Let $S$ be a semigroup. Suppose  $a,b\in S$ are group elements.
 Then $a\sim b$ implies $a\D b$.
\end{corollary}

\begin{corollary} Let $S$ be a completely regular semigroup. Then
the relations $\sim_p$ and $\sim$ on $S$ coincide. In particular,
$\sim_p$ is an equivalence relation.
\end{corollary}
\begin{corollary} Let $S$ be a band. Then $a\sim b$ if and only if
$a\D b$.
\end{corollary}

\section{The general case}\label{s4}

To prove the results of this Section we will use the results of
the previous Section and one important observation, from which we
start.

\begin{proposition}\label{key}
Let $S$ be a regular semigroup and $a\in S$ a group-bound element.
Then $a\sim ae_a$.
\end{proposition}
\begin{proof}
Let $t$ be the height of $a$ (see the definition in the
preliminaries). For the case $t=1$ the statement is obvious as
$ae_a=a$. Suppose $t\geq 2$. For each $i$, $1\leq i\leq t-1$,
denote by $\alpha_i$ any element which is inverse of $a^i$, and by
$\alpha_t$ the element $(a^t)^{-1}$, which is inverse of $a^t$ in
the group $H_{a^t}$. Put $c_0=a$, $c_i=a^{i+1}\alpha_i$, $1\leq
i\leq t$.  Note that for $0\leq i\leq t-1$ we have
\begin{equation}\label{a1}
a^i\alpha_i\cdot a^{i+1}\alpha_{i+1}=a^i\alpha_ia^i\cdot
a\alpha_{i+1}=a^{i+1}\alpha_{i+1}.
\end{equation}
Let $s=c_i$, $t=a^{i+1}\alpha_{i+1}$. Then using~(\ref{a1}) we
obtain
$$st=a\cdot a^i\alpha_i\cdot a^{i+1}\alpha_{i+1}=a\cdot
a^{i+1}\alpha_{i+1}=c_{i+1};$$ $$ts=a^{i+1}\alpha_{i+1}\cdot
a^{i+1}\alpha_i=a\cdot a^i\alpha_i=c_i.$$ It follows that
$c_i\sim_p c_{i+1}$, $0\leq i\leq t-1$. Therefore, $a=c_0\sim
c_t=ae_a$.
\end{proof}
\begin{corollary}\label{c5}
Let the semigroup $S$ be regular and $a,b\in S$ be group-bound
elements. Then $a\sim b$ if and only if $ae_a\sim be_b$.
\end{corollary}

\begin{theorem}\label{main}
Let $S$ be a regular epigroup and $a,b\in S$. Then $a\sim b$ if
and only if there exists a pair of mutually inverse elements
$u,v\in S$ such that $ae_a=u\cdot be_b\cdot v$ and $be_b=v\cdot
ae_a\cdot u$.
\end{theorem}
\begin{proof} The statement follows from Corollary~\ref{c1},
Theorem~\ref{pp1} and Corollary~\ref{c5}.
\end{proof}

\begin{corollary} Let $S$ be a regular semigroup with the zero
element $0$. Then any two nilpotent elements are conjugate, and if
$a\sim b$ and $a$ is nilpotent then $b$ is also nilpotent.
\end{corollary}
\begin{proof} Suppose that $a,b$ are nilpotent. Then $ae_a=be_b=0$,
so that $a\sim b$ by Theorem~\ref{main}.

Suppose now that $a$ is nilpotent and $a\sim b$. This and
Corollary~\ref{c5} imply that $be_b\sim ae_a=0$. It follows now
from Corollary~\ref{c2} that $be_b\D 0$. Thus $be_b=0$, and hence
$e_b\H be_b=0$, so that $e_b=0$. Therefore, $b$ is nilpotent.
\end{proof}

Recall (see~\cite{How}, p.199) that an inverse semigroup $S$ with
the group of units $G$ is called {\em factorizable} provided that
for each $s\in S$ there is $g\in G$ such that $s\leq g$ with
respect to the natural partial order on $S$ i.e.
$ss^{-1}=sg^{-1}$.

The following theorem provides a characterization of conjugacy in
terms of $G$-conjugacy for the class of factorizable inverse
epigroups.

\begin{theorem}\label{p3} Let $S$ be a factorizable inverse epigroup with the identity element $e$ and
the group of units $G$. Let $a,b\in S$. Then $a\sim b$ if and only
if $ae_a\sim_G be_b$.
\end{theorem}
\begin{proof} Since $\sim_G\subset \sim$ and in view of
Corollary~\ref{c5} it is enough to prove only that $a\sim b$
implies $ae_a\sim_G be_b$. Suppose $a\sim b$. It follows from
Theorem~\ref{main} that $ae_a=sbe_bt$ and $be_b=tae_as$ for some
mutually inverse $s\in R_a\cap L_b$ and $t\in R_b\cap L_a$. Since
$S$ is an inverse semigroup it follows that $t$ coincides with
$s^{-1}$
--- the (unique) element, inverse to $s$. That $s^{-1}s\H b$ yields
$s^{-1}sbs^{-1}s=b$. Let $g\in G$ be such that $s\leq g$. Then
$a=sbs^{-1}=ss^{-1}sbs^{-1}ss^{-1}=gs^{-1}sbs^{-1}sg^{-1}=gbg^{-1}$,
and the proof is complete.
\end{proof}

\section{Some examples}\label{appl}

\subsection{Finite transformation semigroups $\IS_n$, $\T_n$ and $\PT_n$}

Let $\IS_n$ be the {\em full finite inverse symmetric semigroup},
i.e. the semigroup of all partial permutations over an $n$-element
set $X=\{1,\dots, n\}$.  The group of units  of $\IS_n$ is the
full symmetric group $\S_n$of all everywhere defined permutations.
The idempotents of $\IS_n$ are precisely the identity maps on
subsets of $X$. Let $\pi\in \IS_n$. Set $G_{\pi}$ to be the
directed graph whose set of vertices $V(G_{\pi})$ coincides with
$X$, and $(x,y)\in E(G_{\pi})$ if and only if $\pi(x)=y$. The
graph $G_{\pi}$ is called the {\em graph of action} of $\pi$.
There are two types of connected components of $G_{\pi}$: {\em
cycles and chains} (see~\cite{Li, GK, KM}). The {\em cyclic type}
and the {\em chain type} of $\pi$ are respectively the (unordered)
tuples $(n_1,\dots, n_k)$ and $(l_1, \dots, l_t)$, where
$n_1,\dots, n_k$ are the lengthes of the cycles of $\pi$, and
$l_1,\dots, l_t$ are the lengthes of the chains of $\pi$ (by the
length of a chain $a_1\to a_2\to\dots \to a_n\to \varnothing$ we
mean here the number $n$ of its vertices). The following Lemma is
straightforward.
\begin{lemma}\label{lll}
\begin{enumerate}
\item $\pi\in \IS_n$ is a group element if and only if all its
chains are trivial, i.e. of length $1$.

\item If the cyclic and chain types of $\pi$ are respectively $(n_1,\dots, n_k)$ and $(l_1, \dots,
l_t)$ then the cyclic and chain types of $\pi e_{\pi}$ are
respectively $(n_1,$ $\dots, n_k)$ and $(1, \dots, 1)$.

\item $\pi\sim_{S_n} \tau$ if and only if the graphs $G_{\pi}$ and
$G_{\tau}$ are isomorphic as directed graphs, which is the case if
and only if the cyclic and chain types of $\pi$ and $\tau$
coincide.
\end{enumerate}
\end{lemma}

Since $\IS_n$ is factorizable, Theorem~\ref{p3} is applicable.
Together with Lemma~\ref{lll} it gives the following criterion of
conjugacy for the semigroup $\IS_n$.

\begin{theorem}[~\cite{Li, GK}]\label{gk}
Let $\pi$, $\tau\in \IS_n$. Then $\pi\sim \tau$ if and only if the
cyclic types of $\pi$ and $\tau$ coincide.
\end{theorem}

Let $\T_n$ and $\PT_n$ be the the semigroups of respectively all
transformations and of all partial transformations (in both cases
not necessarily injective) of the set $X=\{1,\dots, n\}$. Both of
these semigroups are regular while not inverse. In the same vein
as it was done in the case of $\IS_n$ we define the graph of
action $G_{\pi}$ for $\pi\in \T_n$ or $\pi\in\PT_n$. Let $\pi\in
\T_n$ or $\pi\in \PT_n$. One can easily make sure that each
connected component of $G_{\pi}$ contains no more than one cycle.
By the cyclic type of $\pi$ we will mean the (unordered) tuple
$(n_1,\dots, n_k)$, where $n_1\dots, n_k$ are the lengthes of
cycles of $\pi$. Denote the range of $\pi$ by $\ran \pi$, and the
kernel of $\pi$ by $\ker \pi$. Recall that the kernel of $\pi$ is
such a partition of the domain of $\pi$ that $a,b\in X$ belong to
the same block if and only if $a\pi =b\pi $.

\begin{lemma}\label{lema} Let $\pi \in \T_n$ or $\pi\in \PT_n$.
Then
\begin{enumerate}
\item The cyclic types of $\pi$ and $\pi e_{\pi}$ coincide.
\item $\pi$ is a group element if and only if $\ran \pi$ is a
transversal of $\ker \pi$. In the latest case the restriction
${\overline \pi}$ of $\pi$ to $\ran \pi$ is a permutation on the
set $\ran \pi$ and is a group element of $\IS_n$.
\item If $\pi$ is a group element then the cyclic types of $\pi$
and ${\overline \pi}$ coincide.
\item Two group elements $\pi$, $\tau$ $\in \T_n$ (or $\PT_n$) are conjugate if
and only if ${\overline \pi}$ and ${\overline \tau}$ are conjugate
in $\IS_n$.
\end{enumerate}
\end{lemma}

\begin{proof}
{\it 1}. Let ${\mathrm {stran}}\pi=\cap_{k\geq 1} \ran \pi^k$ be
the stable range of $\pi$. For $a\in X$ we have that $a\in
{\mathrm {stran}}\pi$ if and only if $a$ belongs to a cycle in
$G_{\pi}$. This and that $e_{\pi}$ acts identically on ${\mathrm
{stran}}\pi$ imply that $\pi$ and $\pi e_{\pi}$ have the same
cycles (see also~\cite{KM}).

{\it 2}. Follows from the description of Green's relations in
$\T_n$ and $\PT_n$ (see, for example,~\cite{Hig}).

{\it 3}. Since $\ran \pi={\mathrm {stran}}\pi$ in the case of the
group element $\pi$, it follows that the graph of action
$G_{{\overline \pi}}$ of ${\overline \pi}$ is the union of cycles
of $G_{\pi}$, whence the cyclic types of $\pi$ and ${\overline
\pi}$ coincide.

{\it 4}. Let first $\pi\sim \tau$. It follows from
Theorem~\ref{main} and its proof that there are mutually inverse
elements $t\in R_{\pi}\cap L_{\tau}$ and $t'\in R_{\tau}\cap
R_{\pi}$ such that $\pi=t\tau t'$ and $\tau=t'\pi t$. It follows
from the description of Green's relations on $\T_n$ and $\PT_n$
that
$$\ker\,t=\ker \tau,\,\,\, \ker \,t'=\ker \pi, \,\,\, \ran \,t=\ran \pi, \,\,\,
\ran\, t' = \ran \tau.$$

Let ${\overline t}$ be the restriction of $t$ to $\ran\tau$ and
${\overline t'}$ --- the restriction of $t'$ to $\ran \pi$. It
follows that ${\overline t}{\overline t'}$ is the identity map on
$\ran\pi$ and ${\overline t'}{\overline t}$ is the identity map on
$\ran\tau$. Hence, ${\overline t}$ and ${\overline t'}$ are the
pair of mutually inverse elements from $\IS_n$. This, ${\overline
\pi}={\overline t}{\overline \tau}{\overline t'}$, ${\overline
\tau}={\overline t'}{\overline \pi}{\overline t}$ and
Theorem~\ref{main} imply that ${\overline \pi}$ and ${\overline
\tau}$ are $\IS_n$-conjugate.

Now let  the partial permutations ${\overline \pi}$ and
${\overline \tau}$ be $\IS_n$-conjugate. Then there exist
${\overline t}\in R_{{\overline \pi}}\cap L_{{\overline \tau}}$
and ${\overline t'}\in R_{{\overline \tau}}\cap R_{{\overline
\pi}}$ such that ${\overline \pi}={\overline t} {\overline \tau}
{\overline t'}$ and ${\overline \tau}={\overline t'} {\overline
\pi} {\overline  t}$. Define the elements $t,t'\in \T_n$ ($\PT_n$)
 as follows. Set $t$ to be such that $\ker t=\ker \tau$, $\ran t=
\ran \pi =\ran {\overline t}$ and the restriction of $t$ to $\ran
\tau$ coincides with ${\overline t}$. Similarly, set $t'$ to be
such that $\ker t'=\ker \pi$, $\ran t'= \ran \tau =\ran {\overline
t'}$ and the restriction of $t'$ to $\ran \pi$ coincides with
${\overline t}$. Note that it follows from the definitions of
${\overline t}$ and ${\overline t'}$ that $t$ and $t'$ can be
constructed uniquely and happen to be mutually inverse. Moreover,
the construction of $t$ and $t'$ implies $\pi=t\tau t'$ and $\tau
=t'\pi t$. Then by Theorem~\ref{main} $\pi$ and $\tau$ are
conjugate.
\end{proof}

As a corollary we obtain the criterion of $\T_n$- (or $\PT_n$-)
conjugacy in terms of cyclic types of elements.

\begin{theorem}[~\cite{KM}]
Let $\pi, \tau \in \T_n$ ($\PT_n$). Then $\pi$ and $\tau$ are
$\T_n$- ($\PT_n$-) conjugate if and only if their cyclic types
coincide.
\end{theorem}

\subsection{Full semigroups of linear transformations of a
finitely dimensional vector space}

Let $F$ be a field and $V_n$ be an $n$-dimensional vector space
over $F$. An isomorphism $\varphi:U\to W$, where $U, W$ are some
subspaces of $V_n$ is called a {\em partial automorphism} of $V_n$
with the domain $\dom \varphi=U$ and the range $\ran\varphi=W$.
The set of partial automorphisms of $V_n$ with respect to the
composition of partial automorphisms is an inverse semigroup and
denoted by $\PAut(V_n)$. Let $\varphi\in \PAut(V_n)$. For each
positive integer $k$  we have the inclusions
$$
\dom\varphi\supset \dom\varphi^k\supset \dom \varphi^{k+1}\supset
\{0\},
$$
implying that
$$
n\geq \dim U \geq \dim(\dom\varphi^2)\geq\dots
\geq\dim(\dom\varphi^{k})\geq \dots\geq0.
$$
Since at most $n$ of this inequalities are strict we can assert
that starting from some power $t$ we have $\dom\varphi^{t}=
\dom\varphi^{t+i}$ for each $i\geq 0$. It follows that
$\dom\varphi^{t}=\ran\varphi^t$, so that $\varphi^t\in {\mathrm
{GL}}(\dom\varphi^t)$ is a group element of $\PAut(V_n)$, which
shows that $\PAut(V_n)$ is an epigroup (the notation ${\mathrm
{GL}}(W)$ stands for the full linear group over the subspace $W$).
It is easily proved that $\PAut(V_n)$ is factorizable. From
Theorem~\ref{p3} we derive the following criterion of
$\PAut(V_n)$-conjugacy.

\begin{theorem}\label{paut}
Let $\varphi,\psi\in\PAut(V_n)$. Then $\varphi$ and $\psi$ are
$\PAut(V_n)$-conjugate if and only if $\varphi e_{\varphi}$ and
$\psi e_{\psi}$ are ${\mathrm {GL}}(V_n)$-conjugate.
\end{theorem}

Now switch to the regular semigroups $\End(V_n)$ and $\PEnd(V_n)$
of respectively all endomorphisms and all partial endomorphisms of
$V_n$. The same arguments as in the case of $\PAut(V_n)$ show that
both $\End(V_n)$ and $\PEnd(V_n)$ are epigroups.

\begin{lemma} Let $S$ denote one of the semigroups $\End(V_n)$ or
$\PEnd(V_n)$.
\begin{enumerate}
\item $\pi\in S$ is a group element if and only if
$\dom\pi$ decomposes into the direct sum $\dom\pi=\ran\pi\oplus
\ker\pi$. In the latest case the restriction ${\overline\pi}$ of
$\pi$ to $\ran\pi$ is an automorphism of $\ran\pi$ which is a
group element of $\PAut(V_n)$.

\item Two group elements $\pi$, $\tau$ $\in S$  are $S$-conjugate
if and only if ${\overline \pi}$ and ${\overline \tau}$ are
$\PAut(V_n)$- conjugate.
\end{enumerate}
\end{lemma}
\begin{proof}
{\em 1.} Recall that $\pi\in S$ is a group element if and only if
$H_{\pi}$ contains an idempotent, i.e. some projection map
$e=e(V_1,V_2)$, such that $\dom \,e=V_1\oplus V_2$ and $e$ is a
projecting of $\dom \,e$ onto $V_1$ parallelly to $V_2$. The
statement now follows from the fact that $\pi\H e$ if and only if
$\ker\pi=\ker e$ and $\ran\pi =\ran e$.

{\em 2.} The proof is similar to the proof of the fourth statement
of Lemma~\ref{lema}.
\end{proof}

As a corollary we obtain the criterion of conjugacy (where
$S=\End(V_n)$ or $S=\PEnd(V_n)$) in terms of $G$-conjugacy.

\begin{theorem}[~\cite{KM1} for the case of $\End(V_n)$]
Let $S$ denote one of the semigroups $\End(V_n)$ or $\PEnd(V_n)$
and $\varphi,\psi\in S$. Then $\varphi$ and $\psi$ are
$S$-conjugate if and only if ${\overline{\varphi_e{\varphi}}}$ and
${\overline{\psi e_{\psi}}}$ are ${\mathrm {GL}}(V_n)$-conjugate.
\end{theorem}

\subsection{Partial automatic permutations over a finite alphabet}

Recall that a {\em Mealy automaton over a finite alphabet $X$} is
a triple ${\mathcal A}=(X, \varphi, \psi)$, where $Q$ is the set
of {\em internal states} of the automaton, $\varphi: Q\times X\to
Q$--- its {\em transition function} and $\psi: Q\times X\to X$ --
its {\em output function}. In the case when the functions
$\varphi$ and $\psi$ are everywhere defined the automaton
${\mathcal A}$ is called {\em full}, otherwise it is called {\em
partial}. An automaton ${\mathcal A}$ is called {\em initial} if a
state $q_o\in Q$ is marked as an {\em initial state}. Each initial
automaton $({\mathcal A}, q_0)$ over $X$ defines a (partial)
transformation of the set $X^*$ of all words over $X$ by extending
functions $\varphi$ and $\psi$ to the set $Q\times X^*$ as
follows:
$$
\varphi(q,e)=q,\,\,\,\,\,\varphi(q, wx)=\varphi(\varphi(q,w),x);
$$
$$
\psi(q,e)=e,\,\,\,\,\,\psi(q, wx)=\psi(\varphi(q,w),x),
$$
where $e$ denotes the empty word. Now define the transformation
$f_{{\mathcal A},q_0}:X^*\to X^*$ via
\begin{equation}\label{aa}
f_{{\mathcal
A},q_0}(u)=\psi(q_0,x_1)\psi(\varphi(q_0,x_1),x_2)\psi(\varphi(q_0,x_1x_2),x_3)...,
\end{equation}
where $u=x_1x_2x_3\dots\in X^*$. The expression in the right-hand
side of~(\ref{aa}) is undefined if and only if at least one of the
values of $\varphi$ or $\psi$ in it is undefined. Partial
injective transformations which is defined by some partial initial
automaton is called a {\em partial automatic permutation} (or, in
other terminology, a {\em letter-to-letter transduction}). The set
of all partial automatic permutations over $X$ with respect the
composition of maps is an inverse semigroup which will be denoted
by $\ISA(X)$. Note that a partial automatic permutation is a group
element of $\ISA(X)$ if and only if its graph of action has no
chains of length greater than one. The following Lemma is
straightforward.

\begin{lemma}\label{gb}
An element $f\in\ISA(X)$ is group-bound if and only if the
lengthes of its chains are uniformly bounded.
\end{lemma}

A partial automatic permutation $g\in \ISA(X)$ is said to be {\em
finitary} if there exists $l\geq 0$ such that for every word
$x_1x_2\dots \in X^*$ belonging to the domain of $g$ and its image
$y_1y_2\dots =(x_1x_2\dots)^g$ one has $x_i=y_i$ for all $i\geq
l$. The set $Fin\ISA(X)$ of all finitary partial automatic
permutations is an inverse subsemigroup of $\ISA(X)$ and by
Lemma~\ref{gb} it is an epigroup. The group of units of
$Fin\ISA(X)$ coincides with the group $Fin{\mathcal{SA}}(X)$
consisting of all everywhere defined elements of $Fin\ISA(X)$. It
is easily seen that the semigroup $Fin{\mathcal{SA}}(X)$ is
factorizable. From Theorem~\ref{p3} we obtain the following
conjugacy criterion for this semigroup.
\begin{theorem}\label{last}
Two elements $f,g\in Fin\ISA(X)$ are conjugate with respect to
$\sim$ if and only if $fe_f$ and $ge_g$ are
$Fin{\mathcal{SA}}(X)$-conjugate.
\end{theorem}

Theorem~\ref{main} also gives us criteria of conjugacy for the
regular epigroups $Fin{\mathcal A}(X)$ of all (not necessarily
injective) {\em finitary automatic transformations of $X^*$} and
$Fin{\mathcal{PA}}(X)$ of all {\em partial finitary automatic
transformations of $X^*$}. It easily seen that a statement similar
to Lemma~\ref{lema} holds for these semigroups. Therefore we
obtain the following conjugacy criterion.

\begin{theorem}\label{last1}
Let $S$ denote one of the semigroups $Fin{\mathcal A}(X)$ or
$Fin{\mathcal{PA}}(X)$. Two elements $f,g\in S$ are $S$-conjugate
 if and only if ${\overline{fe_f}}$ and ${\overline{ge_g}}$ are
$Fin{\mathcal{SA}}(X)$-conjugate.
\end{theorem}

\section*{Appendix A}\label{ab}

Here we are going to show that for $|X|\geq 2$ the semigroup
$\ISA(X)$ is not an epigroup. For this we give an example of an
automaton $({\mathcal A}, q_0)$ with four states over a two-letter
alphabet $X=\{0,1\}$ such that $f_{{\mathcal A},q_0}$ is not a
group-bound element of $\ISA(X)$. That finite automata $({\mathcal
A}, q_0)$ such that $f_{{\mathcal A},q_0}$ is not group-bound
exist is rather evident. However, so far as to our knowledge, this
fact has never been indicated in the literature. Besides, from
Corollary~\ref{epi} it follows that a non group-bound element can
not be conjugate to a group element. This and the existence of non
group-bound elements assure that the conjugacy criterion for the
semigroup $\ISA(X)$ announced in~\cite{OS} is incorrect.

\begin{figure}
\centering \includegraphics[width=4.5in, height=5.5in]{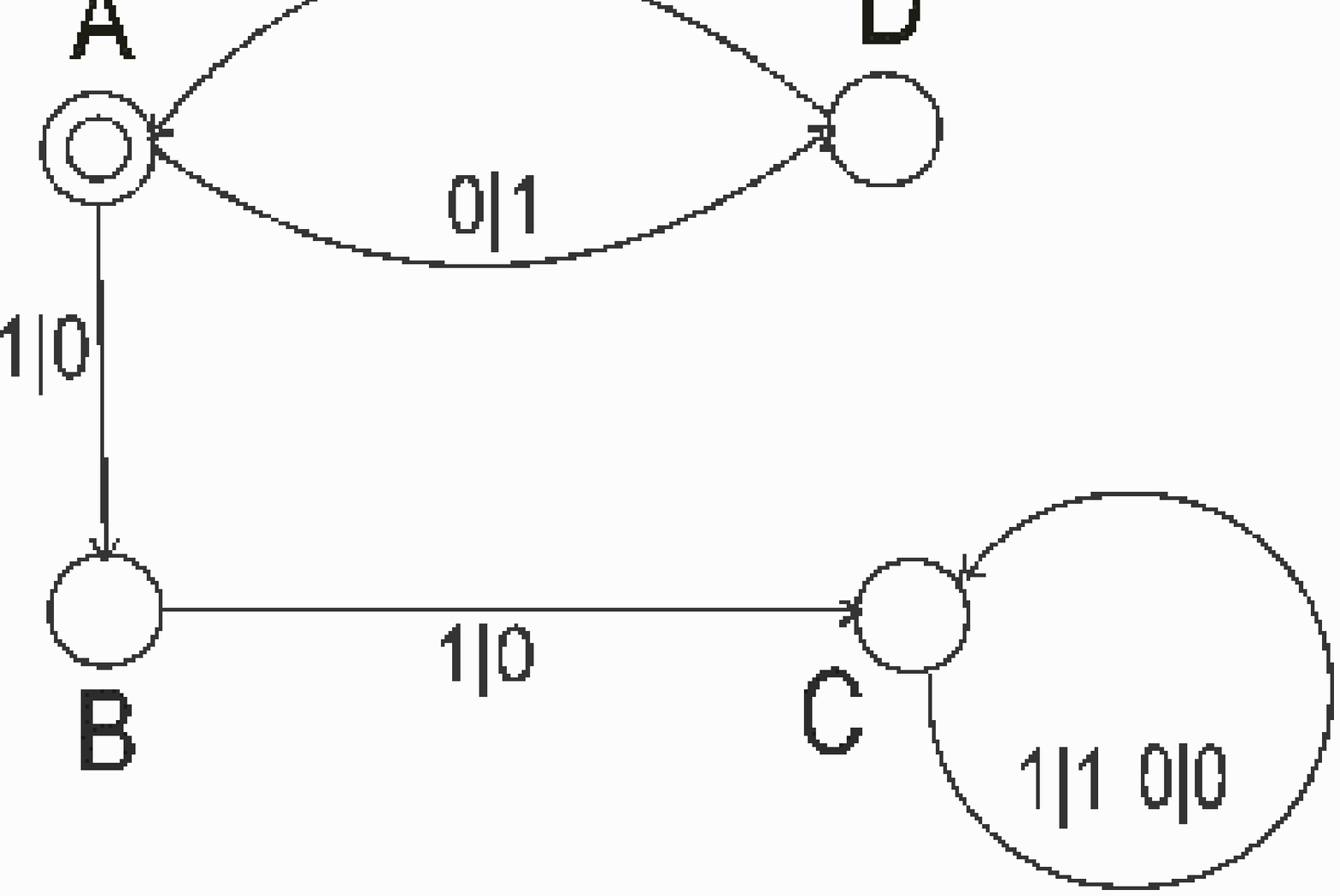}
\vspace*{8pt} \vspace{-7cm} \caption{}\label{fig}
\end{figure}

Construct $({\mathcal A}, q_0)$ $=(Q,\varphi,\psi, q_0)$ as
follows. Let $Q=\{A,B,C,D\}$, $q_0=A$ and
$$
\varphi(A,0)=D,\,\,\,\,\, \varphi(A,1)=B\,\,\,\,\, \ \varphi(B,0)
\text{ undefined},\,\,\,\,\, \varphi(B,1)=C,
$$
$$
\varphi(C,0)=C,\,\,\,\,\,\varphi(C,1)=C,\,\,\,\,\,\varphi(D,0)=A,\,\,\,\,\,\varphi(D,1)
\text{ undefined};
$$
$$
\psi(A,0)=1,\,\,\,\,\, \psi(A,1)=0,\,\,\,\,\, \ \psi(B,0) \text{
undefined},\,\,\,\,\, \psi(B,1)=0,
$$
$$
\psi(C,0)=0,\,\,\,\,\,\psi(C,1)=1,\,\,\,\,\,\psi(D,0)=1,\,\,\,\,\,\psi(D,1)
\text{ undefined}.
$$
The Moore diagram of the constructed automaton is given in
Figure~\ref{fig}, where the initial state is marked by a double
circle, and there is no arrow with the first label $x\in X$
beginning in a state $q\in Q$ if and only if $\varphi(q,x)$ and
$\psi(q,x)$ are undefined.

\begin{lemma}\label{isa} $f_{{\mathcal A},q_0}$ is not a group-bound element of
$\ISA(\{0,1\})$.
\end{lemma}

\begin{proof}
 Show first that for any $k\geq 1$ the orbit (with respect to the
action of $f_{{\mathcal A},q_0}$) of the word ${\underbrace{1\dots
1}_{2k}}$ is a cycle of length $2^k$. Apply induction on $k$. For
$k=1$ we have $11 \mapsto 00 \mapsto 11$. Suppose that $k\geq 1$
and
$$
{\underbrace{1\dots 1}_{2k}}=u_1\mapsto u_2\mapsto \dots \mapsto
u_{2^k}={\underbrace{0\dots 0}_{2k}}\mapsto u_1.
$$
It follows from the definition of $({\mathcal A}, q_0)$ that
\begin{multline*}
{\underbrace{1\dots 1}_{2k+2}}=u_111\mapsto u_211\mapsto \dots
\mapsto u_{2^k}11={\underbrace{0\dots 0}_{2k}}11\mapsto\\
{\underbrace{1\dots 1}_{2k}}00=u_100\mapsto u_200 \mapsto\dots
\mapsto u_{2^k}00={\underbrace{1\dots 1}_{2k+2}}\mapsto u_111,
\end{multline*}
as required. Let now $k\geq 1$. Then
$$
{\underbrace{1\dots 1}_{2k}}01=u_101\mapsto u_201\mapsto \dots
\mapsto u_{2^k}01={\underbrace{0\dots 0}_{2^k}}01,
$$
and $f_{{\mathcal A},q_0}({\underbrace{0\dots 0}_{2k}}01)$ is
undefined. Therefore, the word ${\underbrace{1\dots 1}_{2k}}01$
belongs to the chain of length at least $2^k$, $k\geq 1$. The
statement now follows from Lemma~\ref{gb}.
\end{proof}

\begin{corollary} If $|X|\geq 2$ then in $\ISA(X)$ there are
conjugacy classes without group elements.
\end{corollary}
\begin{proof}
This follows from Corollary~\ref{epi} and Lemma~\ref{isa}.
\end{proof}


\end{document}